%% file: PaperArXiv.tex
\documentclass[11pt]{article}
\usepackage[english]{babel}    
\usepackage[utf8]{inputenc}    
\usepackage{blindtext}         
\usepackage[T1]{fontenc}       

\usepackage[sumlimits,%
                intlimits,%
                namelimits%
                ]{amsmath}         
\usepackage{amsthm}            
\usepackage{amscd}             

\usepackage{fouridx}
\usepackage{mathtools}
\usepackage{tikz}
\usepackage{amssymb}

\usepackage[hyperfootnotes=false]{hyperref}
\usepackage{cleveref}
\usepackage{autonum}
\usepackage{enumerate}

\usepackage{pgfplots}
\pgfplotsset{compat=newest}
\usetikzlibrary{plotmarks}
\usetikzlibrary{arrows.meta}
\usepgfplotslibrary{patchplots}
\usetikzlibrary{math}
\usepackage{tkz-euclide}

\newtheorem{theorem}{Theorem}

\def\RR{{\mathbb{R}}}

\def\phi{\varphi}

  \newcommand{\RN}[1]{%
  \textup{\uppercase\expandafter{\romannumeral#1}}%
}

\begin{document}
\begin{center}
 \bf On Quasi-Interpolation and their associated shift-invariant space using a new class of generalized Thin Plate Splines and Inverse Multiquadrics
 
\end{center}
\bigskip

\begin{center}
  \it Mathis Ortmann, Justus-Liebig University,
  Mathematic Department, 35392 Giessen, and \\ \smallskip
  Martin Buhmann, Justus-Liebig University, Mathematic Department, 35392 Giessen

\end{center}
\bigskip\bigskip
{\small {\it Abstract: \/} A new generalization of  shifted thin plate splines
 $$\varphi(x)=(c^{2d}+||x||^{2d})\log\left(c^{2d}+||x||^{2d}\right),\qquad x\in\mathbb{R}^n, d\in \mathbb{N}, c>0$$ is presented to increase the accuracy of quasi-interpolation further. With the restriction to Euclidean spaces of even dimensionality, the generalization can be used to generate a quasi-Lagrange operator that reproduces all polynomials of degree $n+2d-1$. It thus complements the case of the newly proposed generalized multiquadric $\varphi(x)=\sqrt{c^{2d}+||x||^{2d}},\quad x\in\mathbb{R}^n, d\in \mathbb{N}, c>0$, which is restricted to odd dimensions \cite{ortmann}. 
This generalization improves the approximation order by a factor of $\mathcal{O}\left(h^{2(d-1)}\right)$, where $d=1$ represents the classical thin plate spline. The results are then compared with the theoretical optimal approximation from the shift-invariant space generated by these functions.
Moreover, we introduce a new class of inverse multiquadrics
$$\varphi(x)=\left(c^\lambda +||x||^\lambda\right)^\beta,\qquad x\in\mathbb{R}^n, \lambda \in\mathbb{R},\beta \in \mathbb{R}\backslash\mathbb{N}, c>0. $$
We provide an explicit representation of the generalized Fourier transform and discuss its asymptotic behaviour near the origin. Particular emphasis is placed on the case where $\lambda$ and $\beta$ are both negative. It is demonstrated that, in dimensions $n\geq3$, it is possible to build a quasi-Lagrange operator that reproduces all polynomials of degree $n-3$ when $n$ is even and of degree $\frac{n-1}{2}$ when n is odd. Furthermore, the uniform approximation error is given by $\mathcal{O}\left(h^{n-2}\log(1/h)\right)$ for $n$ even and $\mathcal{O}\left(h^{\frac{n-3}{2}}\right)$ for $n$ odd. Here, $h>0$ denotes the fill distance.

\date{\today}  
 \bigskip\bigskip\bigskip
 \title{Introduction}
\section{Introduction}

Quasi-interpolation is a frequently employed technique in the fields of image processing, surface reconstruction, and data fitting \cite{WOS:000305709900003},\cite{WOS:000513922200014}. It is also utilised in the solution of partial differential equations \cite{CHEN}. The method is based on the theory of shift-invariant spaces, similar to other impactful methods such as wavelets \cite{WaveletTour}, shearlets \cite{Kutyniok.2012} and radial basis functions (RBF) \cite{rbf}. There are several researchers searching for different shift-invariant spaces and compare the convergence orders of the best approximant to given data points generated by an unknown function \cite{HOLTZ200597}, \cite{Ron.1992}. The best approximation from an invariant space may be excessively complex when attempting to compute it. In contrast, quasi-interpolation provides a methodology for the construction of an approximant. \cite{qi}. The convergence rate of the quasi-interpolant typically differs from the best approximant by only a factor $\log \left(1/h\right)$, where $h$ denotes the fill distance.

This paper will look at quasi-interpolants built from a principal shift-invariant space  $\mathcal{S}=\{\phi(\cdot-\beta)| \beta \in \mathbb{Z}^n \}$  generated by a single function, where $\phi$ is a RBF. Classical choices for RBFs are the multiquadric $\phi(x)=\sqrt{c^2+||x||^2}, c>0$, that is used in odd dimensions and the shifted thin plate spline $\phi(||x||)= (c^2+||x||^2)\log( c^2+||x||^2)$ for even dimensions. 
The choice of the RBF has a significant impact on the approximation order of the quasi-interpolant. Consequently, numerous researchers are engaged in the pursuit of novel functions and generalizations. This paper introduces a novel generalization of the shifted thin plate spline
 \begin{equation}
\phi(x)= (c^{2d}+||x||^{2d})\log( c^{2d}+||x||^{2d}),\qquad d\in\mathbb{N}, x\in \mathbb{R}^n, c>0. \label {eq: TPS}
\end{equation}
This generalization complements the recently proposed generalized multiquadric function $\phi(x)=\sqrt{c^{2d}+||x||^{2d}}$, that is presented in \cite{ortmann}.
Furthermore we introduce a new class of inverse multiquadrics 
$$\phi(x)=\left(c^\lambda +||x||^\lambda\right)^\beta,\qquad x\in\mathbb{R}^n, \lambda \in\mathbb{R},\beta \in \mathbb{R}\backslash\mathbb{N}, c>0, $$
which constitutes a further generalization to the multiquadrics described in \cite{ortmann} for $\lambda, \beta >0$. Nevertheless, we discover all possible choices for $\lambda$ and $\beta$. Particular emphasis is given to the case where $\lambda$ and $\beta$ are both negative, which we refer to as the new class of inverse multiquadrics.
The quasi-Lagrange function $\Psi$ is then given by a linear combination of RBFs
\begin{equation}
\Psi(x)=\sum_{\alpha \in \mathbb{Z}^n}\mu_\alpha \phi(x-\alpha),\qquad x\in\mathbb{R}^n,
\end{equation}
where the coeffients $\mu_\alpha$ may be either finite or infinite in number, depending on the singularity of the Fourier transform of the RBF at the origin. The quasi-interpolant is then given by
\begin{equation}
Q_hf(x)=\sum_{j\in\mathbb{Z}^n}f(jh)\Psi(x/h-j),\qquad x\in\mathbb{R}^n.
\end{equation}
The convergence rate for $h\to 0$ can be identified using the well-known theorem from Strang and Fix cited below.

\begin{theorem}\label{SF}[Strang and Fix conditions]

Let $m$ be a positive integer and $\Psi:\mathbb{R}^n \rightarrow \mathbb{R}$ be a function such that
\begin{enumerate}
\item there exists a non-negative real valued $\ell$ such that, when $\|x\|\rightarrow \infty$, $|\Psi(x)|=\mathcal{O}(\Vert x\Vert^{-n-m-\ell})$ and this implies $\hat\Psi \in C^{m+\ell-1}(\mathbb{R}^n),$
\item $D^{\alpha} \hat{\Psi}(0)=0$, $\forall \alpha \in \mathbb{Z}_+^n$, $1\leq \vert \alpha\vert\leq m$, and $\hat{\Psi}(0)=1$,
\item $ D^{\alpha}\hat{\Psi}(2\pi j)=0,\ \forall j \in \mathbb{Z}^n\setminus\{0\}$ and  $\forall \alpha \in \mathbb{Z}_+^n$ with $ \vert\alpha\vert \leq m$. 
\end{enumerate}

Then the quasi-interpolant
\begin{equation}
Q_hf(x)=\sum_{j\in\mathbb{Z}^n}f(jh)\Psi(x/h-j),\qquad x\in\RR^n,
\end{equation}
is well-defined and exact on the space of polynomials of degree $m$ and the uniform approximation error can be estimated by
$$\Vert Q_hf-f\Vert_{\infty}=\begin{cases}\mathcal{O}(h^{m+\ell}), & \text{when } 0<\ell<1,\\
\mathcal{O}(h^{m+1}\log (1/h)), & \text{when } \ell=1,\\
\mathcal{O}(h^{m+1}), & \text{when } \ell>1,\\
\end{cases}$$
for $h\rightarrow 0$ and a bounded function $f\in C^{m+1}(\RR^n)$ with bounded derivatives \cite{BUHMANN2015156}.
\end{theorem}

\section{Fourier transform}
 This section presents a detailed calculation of the Fourier transform of the novel class of generalized thin plate splines.
To employ the \textit{Strang and Fix conditions}, one needs the Fourier transform the quasi-Lagrange function $\Psi$ and thus in particular the Fourier transform of the RBF. The $n$-dimensional Fourier transform is defined as
\begin{align}
\hat f(y):=\int_{\mathbb{R}^n}f(x)\textnormal{e}^{-\mathrm{i}x\cdot y} \mathrm{d}x ,\qquad y\in\mathbb{R}^n \, ,
\end{align}
where $f$ is an absolutely integrable function. Since $\phi$ in equation~(\ref{eq: TPS}) is not integrable, we need to employ the theory of generalized Fourier transforms \cite{jones}.

\begin{theorem}\label{FT}
The generalized $n$-dimensional Fourier transform of
 \begin{align}
\phi(x)&=(c^{2d}+||x||^{2d})\log\left(c^{2d}+||x||^{2d}\right), \qquad  x \in \mathbb{R}^n,\, r=||x||, \,d\in \mathbb{N},\, c>0\\
&=\phi_\RN{1}(x)+\phi_\RN{2}(x)+\phi_\RN{3}(x)+\phi_\RN{4}(x),
\intertext{with}
\phi_\RN{1}(x)&=2d c^{2d}\log\left(c\right),\\
\phi_\RN{2}(x)&=2d r^{2d}\log\left (c\right),\\
\phi_\RN{3}(x)&=r^{2d}\log\left(1+\left (\frac{r}{c}\right )^{2d}\right),\\
\phi_\RN{4}(x)&=c^{2d}\log\left(1+\left (\frac{r}{c}\right )^{2d}\right)
\end{align}

 is given by
\begin{align}
\hat\phi(y)&=\hat\phi_\RN{1}(y)+\hat\phi_\RN{2}(y)+\hat\phi_\RN{3}(y)+\hat\phi_\RN{4}(y),\\
\intertext{where}
\hat\phi_\RN{1}(y)&=2d c^{2d}\log\left(c\right)\delta (y), \\
\hat\phi_\RN{2}(y)&=\begin{dcases}
2d\log(c)\frac{\Gamma\left(d+n/2\right)}{\Gamma(-d)}2^{2d+n}\pi^{n/2}||y||^{-2d-n},\qquad &d\notin \mathbb{N} \land d\notin -\mathbb{N}-\frac{n}{2}\\
2d\log(c)\left(2\pi\right)^n\left(-1\right)^d\left(\partial_1^2+\cdots +\partial_n^2\right)^d \delta(y),\qquad &d\in \mathbb{N}\\
2d\log(c)\frac{(-1)^{\frac{-2d-n}{2}}\pi^{n/2}||y||^{-2d-n}}{\Gamma(d)\Gamma(\frac{2d-n-2}{2})2^{\frac{-2d-n}{2}}} \times \\ \times\left(\Psi^0\left(d-1\right)+\Psi^0\left(d-n/2\right)+2\log\left(\frac{2}{||y||}\right)			\right), &d\in -\mathbb{N}-\frac{n}{2},
\end{dcases}\\
\hat\phi_\RN{3}(y)&=\frac{2^{n+2d}\pi^{n/2}}{||y||^{n+2d}}H_{2,4}^{\,3,1}\!\left(\left.{\begin{matrix}\left(0,1\right),\left(1,1\right)\\ \left(0,1\right),\left(0,1\right),\left(\frac{n}{2}+d,d\right),\left(1+d,d\right)\end{matrix}}\;\right|\,\left(\frac{c||y||}{2}\right )^{2d}\right),\\
\hat\phi_\RN{4}(y)&=\frac{2^{n}\pi^{n/2}c^{2d}}{||y||^{n}}H_{2,4}^{\,3,1}\!\left(\left.{\begin{matrix}\left(0,1\right),\left(1,1\right)\\ \left(0,1\right),\left(0,1\right),\left(\frac{n}{2},d\right),\left(1,d\right)\end{matrix}}\;\right|\,\left(\frac{c||y||}{2}\right )^{2d}\right)
\end{align}
are the Fourier transforms. 
The definition of the Fox $H$-function is given in the proof.
\end{theorem}

\subsection*{\textit{Notes}}
This coincides for $d=1$ with the classical case
\begin{align}
\hat\phi_\RN{1}(y)&=2 c^{2}\log\left(c\right)\delta (y)\\
\hat\phi_\RN{2}(y)&=-2\log(c)\left(2\pi\right)^n\left(\partial_1^2+\cdots +\partial_n^2\right) \delta(y)\\
\hat\phi_\RN{3}(y)&=2^{n+2} \pi ^{n/2} s^{-n-2} G_{1,3}^{3,0}\left(\frac{c^2 s^2}{4}|
\begin{array}{c}
 1 \\
 0,2,\frac{n+2}{2} \\
\end{array}
\right) \\
\hat\phi_\RN{4}(y)&=c \left(-2^{\frac{n}{2}+1}\right) \pi ^{n/2} \left(\frac{c}{s}\right)^{n/2} K_{\frac{n}{2}}(c s)\\
\hat\phi(y)&=\hat\phi_\RN{3}(y)+\hat\phi_\RN{4}(y)\\
&=4 (2 \pi )^{n/2} \left(\frac{c}{s}\right)^{\frac{n}{2}+1} K_{\frac{n}{2}+1}(c s), \qquad ||y||>0
\intertext{and even further in the limit $ c\to 0, d \in \mathbb{N} $, $\hat\phi$ converges to}
\hat\phi_\RN{1}(y)&=0 ,\qquad & ||y||>0,\\
\hat\phi_\RN{2}(y)&=0, & ||y||>0,\\
\hat\phi_\RN{3}(y)&=(-1)^{d+1} 2^{n+2d} \pi ^{n/2}d\,   \Gamma \left(d+1\right) \Gamma \left(\frac{n}{2}+d\right) ||y||^{-n-2d}, & ||y||>0,\\
\hat\phi_\RN{4}(y)&=0, & ||y||>0,
\intertext{which is in total agreement with the results given in \cite{jones}.}
\end{align}

\subsection*{\textit{Proof}}
The Fourier transform of a radially symmetric function is also a radially symmetric function. Thus we can write

\begin{align}
\hat\phi (s)= \lim_{\varepsilon \to 0_+}\frac{(2\pi)^{\frac{n}{2}}} {s^{\frac{n}{2}-1}}\int_0^\infty r^\frac{n}{2}\phi(r) J_{\frac{n}{2}-1}(sr) \mathrm{e}^{-\varepsilon r^2} \, \mathrm{d}r, \qquad s\in \mathbb{R}_+ ,
\end{align}
where we introduced the Gauss function as a convergence generating factor to keep the integral finite. 
Normally this refers to the Hankel transform of order $\frac{n}{2}-1$, but here it is more appropriate to think of the integral as a Mellin transform.
First we rewrite the RBF into four independent terms
\begin{align}
\phi(r)=\underbrace{2d c^{2d}\log\left(c\right)}_{\phi_\RN{1}(r)}+\underbrace{2d r^{2d}\log\left (c\right)}_{\phi_\RN{2}(r)}+\underbrace{r^{2d}\log\left(1+\left (\frac{r}{c}\right )^{2d}\right)}_{\phi_\RN{3}(r)}+\underbrace{c^{2d}\log\left(1+\left (\frac{r}{c}\right )^{2d}\right)}_{\phi_\RN{4}(r)}
\end{align}
and examine each term independently.
The first term $\phi_\RN{1}(x)$ is just a constant, so the generalized Fourier transform is given by the dirac delta distribution.
The second term can be found in \cite{jones} and is given by
\begin{align}
\hat\phi_\RN{2}(y)&=\begin{dcases}
2d\log(c)\frac{\Gamma\left(d+n/2\right)}{\Gamma(-d)}2^{2d+n}\pi^{n/2}s^{-2d-n},\qquad &d\notin \mathbb{N} \land d\notin -\mathbb{N}-\frac{n}{2}\\
2d\log(c)\left(2\pi\right)^n\left(-1\right)^d\left(\partial_1^2+\cdots +\partial_n^2\right)^d \delta(y),\qquad &d\in \mathbb{N}\\
2d\log(c)\frac{(-1)^{\frac{-2d-n}{2}}\pi^{n/2}s^{-2d-n}}{\Gamma(d)\Gamma(\frac{2d-n-2}{2})2^{\frac{-2d-n}{2}}} \times \\ \times\left(\Psi^0\left(d-1\right)+\Psi^0\left(d-n/2\right)+2\log\left(\frac{2}{s}\right)			\right), &d\in -\mathbb{N}-\frac{n}{2},
\end{dcases}\\
\end{align}
where $s=||y||$ is the radial part of the argument $y$.
The calculation of the Fourier transform of $\phi_\RN{3}$ and $\phi_\RN{4}$ requires more effort. Since the steps are quite similar we will perform them simultaneously. To do so we rewrite the Bessel function and the logarithm into Meijer $G$-functions as follows:

\begin{align}
\log\left(1+\left (\frac{r}{c}\right )^{2d}\right)&=G_{2,2}^{\,1,2}\!\left(\left.{\begin{matrix}1,1\\ 1,0\end{matrix}}\;\right|\,\left(\frac{r}{c}\right )^{2d}\right),\\
J_{\frac{n}{2}-1}(sr)&=G_{0,2}^{\,1,0}\!\left(\left.{\begin{matrix}-\\ \frac{n-2}{4},-\frac{n-2}{4}\end{matrix}}\;\right|\,\frac{(sr)^2}{4}\right), \, r>0.
\end{align}
The Meijer $G$-function is defined as
\begin{align}
G_{p,q}^{\,m,n}\!\left(\left.{\begin{matrix}a_{1},\dots ,a_{p}\\b_{1},\dots ,b_{q}\end{matrix}}\;\right|\,z\right)={\frac {1}{2\pi \mathrm{i}}}\int _{L}{\frac {\prod _{j=1}^{m}\Gamma (b_{j}-s)\prod _{j=1}^{n}\Gamma (1-a_{j}+s)}{\prod _{j=m+1}^{q}\Gamma (1-b_{j}+s)\prod _{j=n+1}^{p}\Gamma (a_{j}-s)}}\,z^{s}\,\mathrm{d}s,
\end{align}
where the path of integration can have three different shapes \cite{handbook}, depending on the parameters $a_1,\dots,a_p$ and $b_1,\dots,b_q$. 
The Meijer $G$-function is a Mellin-Barnes type integral and can be viewed as an inverse Mellin transform \cite{Carlos}, because
\begin{align}
\int _{0}^{\infty }z^{s-1}\;G_{p,q}^{\,m,n}\!\left(\left.{\begin{matrix}a_{1},\dots ,a_{p}\\b_{1},\dots ,b_{q} \end{matrix}}\;\right|\,\eta z\right)\mathrm{d}z={\frac {\eta ^{-s}\prod _{j=1}^{m}\Gamma (b_{j}+s)\prod _{j=1}^{n}\Gamma (1-a_{j}-s)}{\prod _{j=m+1}^{q}\Gamma (1-b_{j}-s)\prod _{j=n+1}^{p}\Gamma (a_{j}+s)}},
\end{align}
where $s,\eta \in\mathbb{R}$.
Substituting $\tilde r =r^2$ and using the identity
\begin{align}
z^{\rho }\;G_{p,q}^{\,m,n}\!\left(\left.{\begin{matrix}\mathbf {a_{p}} \\\mathbf {b_{q}} \end{matrix}}\;\right|\,z\right)=G_{p,q}^{\,m,n}\!\left(\left.{\begin{matrix}\mathbf {a_{p}} +\rho \\\mathbf {b_{q}} +\rho \end{matrix}}\;\right|\,z\right),
\end{align}
where $\mathbf {a_{p}}=[a_1,\dots, a_p]$ and $\mathbf {b_{q}}=[b_1,\dots,b_q]$ are the vectors containing the arguments of the Meijer $G$-function, on the Bessel functions yields

\begin{align}
\hat\phi_\RN{3} (s)&=  \frac{2^{n}\pi^{\frac{n}{2}}} {4s^{n-2}}\lim_{\varepsilon \to 0_+}\int_0^\infty  \mathrm{e}^{-\varepsilon \tilde r}\tilde r^d G_{2,2}^{\,1,2}\!\left(\left.{\begin{matrix}1,1\\ 1,0\end{matrix}}\;\right|\,\frac{\tilde r^d}{c^{2d}}\right) G_{0,2}^{\,1,0}\!\left(\left.{\begin{matrix}-\\ \frac{n}{2}-1,0\end{matrix}}\;\right|\,\frac{s^2\tilde r}{4}\right) \mathrm{d}\tilde r \qquad \textnormal{and}\\
\hat\phi_\RN{4} (s)&=  \frac{2^{n}\pi^{\frac{n}{2}}c^{2d}} {4s^{n-2}}\lim_{\varepsilon \to 0_+}\int_0^\infty  \mathrm{e}^{-\varepsilon \tilde r}  G_{2,2}^{\,1,2}\!\left(\left.{\begin{matrix}1,1\\ 1,0\end{matrix}}\;\right|\,\frac{\tilde r^d}{c^{2d}}\right) G_{0,2}^{\,1,0}\!\left(\left.{\begin{matrix}-\\ \frac{n}{2}-1,0\end{matrix}}\;\right|\,\frac{s^2\tilde r}{4}\right) \mathrm{d}\tilde r.
\end{align}
Writing the first Meijer $G$-function as a Mellin-Barnes integral, changing the order of integration results in
\begin{align}
\hat\phi_\RN{3} (s)&= \frac{2^{n}\pi^{\frac{n}{2}}} {4s^{n-2}}\lim_{\varepsilon \to 0_+} \frac{1}{2\pi \mathrm{i}} \int_L\frac{\Gamma \left(1-t\right)\Gamma \left(t\right)^2}{\Gamma \left(1+t\right)c^{2dt}}  
 \int_0^\infty\tilde r^{d(t+1)}  \mathrm{e}^{-\varepsilon \tilde r}
G_{0,2}^{\,1,0}\!\left(\left.{\begin{matrix}-\\ \frac{n}{2}-1,0\end{matrix}}\;\right|\,\frac{s^2\tilde r}{4}\right) \mathrm{d}\tilde r \, \mathrm{d}t,\\
\hat\phi_\RN{4} (s)&= \frac{2^{n}\pi^{\frac{n}{2}}c^{2d}} {4s^{n-2}}\lim_{\varepsilon \to 0_+} \frac{1}{2\pi \mathrm{i}} \int_L\frac{\Gamma \left(1-t\right)\Gamma \left(t\right)^2}{\Gamma \left(1+t\right)c^{2dt}}  
 \int_0^\infty\tilde r^{dt}  \mathrm{e}^{-\varepsilon \tilde r}
G_{0,2}^{\,1,0}\!\left(\left.{\begin{matrix}-\\ \frac{n}{2}-1,0\end{matrix}}\;\right|\,\frac{s^2\tilde r}{4}\right) \mathrm{d}\tilde r \, \mathrm{d}t.
\end{align}
The change of integration is valid through Fubini's theorem. Additionally, the inner integral is a Laplace transform and can be found in \cite{Prudnikov}, providing
\begin{align}
\hat\phi_\RN{3} (s)&= \frac{2^{n}\pi^{\frac{n}{2}}} {4s^{n-2}}\lim_{\varepsilon \to 0_+} \frac{1}{2\pi \mathrm{i}} \int_L\frac{\Gamma \left(1-t\right)\Gamma \left(t\right)^2}{\Gamma \left(1+t\right)c^{2dt}}  
 \varepsilon^{-dt-d-1} G_{1,2}^{\,1,1}\!\left(\left.{\begin{matrix}-dt-d\\ \frac{n}{2}-1,0\end{matrix}}\;\right|\,\frac{s^2}{4\varepsilon}\right) \mathrm{d}t,\\
\hat\phi_\RN{4} (s)&= \frac{2^{n}\pi^{\frac{n}{2}}c^{2d}} {4s^{n-2}}\lim_{\varepsilon \to 0_+} \frac{1}{2\pi \mathrm{i}} \int_L\frac{\Gamma \left(1-t\right)\Gamma \left(t\right)^2}{\Gamma \left(1+t\right)c^{2dt}}  
 \varepsilon^{-dt-1} G_{1,2}^{\,1,1}\!\left(\left.{\begin{matrix}-dt\\ \frac{n}{2}-1,0\end{matrix}}\;\right|\,\frac{s^2}{4\varepsilon}\right)\mathrm{d}t.\\
\end{align}
Using the Beppo Levi monotone convergence theorem, the integral and the limit can be exchanged. To evaluate the limit, we  write the Meijer $G$-function as a Mellin-Barnes integral and expand it into a series using the residue theorem. It is easy to see that in the limit $\varepsilon \to 0_+$ the only nonvanishing term is induced by the first pole at $dt+1$ for $\hat\phi_\RN{3}$ and at $d(t+1)+1$ for $\hat \phi_\RN{4}$. Ultimately, we obtain the integral representations

\begin{align}
\hat\phi_\RN{3} (s)&= \frac{2^{n+2d}\pi^{\frac{n}{2}}} {s^{n+2d}} \frac{1}{2\pi \mathrm{i}} \int_L\frac{\Gamma \left(1-t\right)\Gamma \left(t\right)^2\Gamma\left(n/2+dt+d\right)}{\Gamma \left(1+t\right)\Gamma\left(-dt-d\right)}\left(\frac{2}{cs}\right)^{2dt}\mathrm{d}t \label{eq: FT3}\\
\hat\phi_\RN{4} (s)&= \frac{2^{n}\pi^{\frac{n}{2}}c^{2d}} {s^{n}} \frac{1}{2\pi \mathrm{i}} \int_L\frac{\Gamma \left(1-t\right)\Gamma \left(t\right)^2\Gamma\left(n/2+dt\right)}{\Gamma \left(1+t\right)\Gamma\left(-dt\right)}\left(\frac{2}{cs}\right)^{2dt}\mathrm{d}t. \label{eq: FT4}
\end{align}

The sum of the integral representations of $\hat\phi_\RN{3} (s)+\hat\phi_\RN{4} (s)$ can be further simplified to
\begin{align}
\hat\phi_\RN{3} (s)+\hat\phi_\RN{4} (s)=\frac{2^{n}\pi^{\frac{n}{2}}c^{2d}} {s^{n}} \frac{1}{2\pi \mathrm{i}} \int_L\frac{\Gamma \left(-t\right)\Gamma \left(t-1\right)\Gamma\left(n/2+dt\right)}{\Gamma\left(-dt\right)}\left(\frac{2}{cs}\right)^{2dt}\mathrm{d}t, \label{eq: FT3+4}
\end{align}
where the path $L$ can be parametrized by $t=r+\mathrm{i}\tilde t$ for some $r\in(1,\infty)$ and $\tilde t \in (-\infty ,\infty)$.
Comparing this with the definition of the Fox $H$-function
\begin{align}
 &H_{p,q}^{\,m,n}\!\left(\left.{\begin{matrix}\left(a_1,A_1\right),\dots,\left(a_p,A_p\right)\\ \left(b_1,B_1\right),\dots, \left(b_q,B_q\right)\end{matrix}}\;\right|\,z\right)\\
&={\frac {1}{2\pi \mathrm{i}}}\int _{L}{\frac { \prod _{j=1}^{m}\Gamma (b_{j}+B_{j}s)\,\prod _{j=1}^{n}\Gamma (1-a_{j}-A_{j}s)}{ \prod _{j=m+1}^{q}\Gamma (1-b_{j}-B_{j}s)\,\prod _{j=n+1}^{p}\Gamma (a_{j}+A_{j}s)}}z^{-s}\,ds,
\end{align}
where $A_1\dots A_p$ and $B_1\dots B_q$ are positive numbers, the desired representation of the Fox $H$-function is obtained.


\section{Asymptotic behaviour}
Next we analyze the asymptotic behaviour of $\hat\phi(s)$ as $s\to 0_+$. 

\begin{theorem}\label{Satz: asymptoticbehaviour}
Let $\hat\phi(s)$ be the generalized Fourier transform of $\phi(x)=(c^{2d}+||x||^{2d})\log\left(c^{2d}+||x||^{2d}\right) $, where $n$ is the dimension and $d \in \mathbb{R}_+$ is the generalization parameter as computed above. Then

\begin{align}
\hat\phi_\RN{3}(s)= 
\begin{cases}
\mathcal{O}\left(s^{-n-2d}\right)\qquad &\textnormal{as } s\to 0 \textnormal{ for } d \in \mathbb{N}, 0<s<1 \\
\mathcal{O}\left(s^{-n-2d}\log\left(s^{-1}\right)\right) \qquad &\textnormal{as } s\to 0 \textnormal{ for } d\notin \mathbb{N}, 0<s<1
\end{cases}
\end{align}
The exact asymptotic behaviour for $s \to 0_+$ can be found in equations ~(\ref{eq: d in N}) and (\ref{eq: d not in N}). Note that the only proposed cases are those where $n+2d$ is an integer.
\end{theorem}

\subsection*{\textit{Proof}}

The asymptotic behaviour of $\hat\phi_\RN{3}$ is determined from equation~(\ref{eq: FT3}). The path $L$ is a loop that starts at infinity on a line parallel to the negative real axis, encircles the poles of the $\Gamma\left(t\right)$ once in the positive sense and returns to infinity on another line parallel to the negative real axis. So the first pole which will determine the asymptotic behaviour is at $t_0=0$. Given this, we can simplify equation~(\ref{eq: FT3}) to
\begin{align}
\hat\phi_\RN{3} (s)&= \frac{2^{n+2d}\pi^{\frac{n}{2}}} {s^{n+2d}} \frac{1}{2\pi \mathrm{i}} \int_L\frac{\pi}{t\sin\left(\pi t\right)}\frac{\Gamma\left(n/2+dt+d\right)}{\Gamma\left(-dt-d\right)}\left(\frac{2}{cs}\right)^{2dt}\mathrm{d}t\, .
\end{align}

\subsubsection*{Case \RN{1}: $d\in \mathbb{N}$ }
If $d$ is an integer the integral has a simple pole at zero. The evaluation using the residue theorem can be easily done and provides
\begin{align}
\hat\phi_\RN{3} (s)\sim (-1)^{d+1} \frac{2^{n+2d}\pi^{\frac{n}{2}}} {s^{n+2d}} d\, \Gamma\left(d+1\right )\Gamma\left(\frac{n}{2}+d\right)\qquad \textnormal{as } s\to 0_+. \label{eq: d in N}
\end{align}

\subsubsection*{Case \RN{2}: $d\notin \mathbb{N}$ }
If $d$ is not an integer, then the pole at zero is of order two. Therefore one has to evaluate
\begin{align}
\lim_{t\to 0}\frac{\partial}{\partial t} \frac{2^{n+2d}\pi^{\frac{n}{2}}} {s^{n+2d}} \frac{t^2\pi}{t\sin\left(\pi t\right)}\frac{\Gamma\left(n/2+dt+d\right)}{\Gamma\left(-dt-d\right)}\left(\frac{2}{cs}\right)^{2dt}.
\end{align}
Straightforward calculations lead to
\begin{align}
\hat\phi_\RN{3} (s)\sim  -\frac{2^{n+2d}\pi^{\frac{n}{2}}} {s^{n+2d}}\frac{d\Gamma\left(d+\frac{n}{2}\right)}{\Gamma\left(-d\right)}\left( 2\log\left(\frac{cs}{2}\right)-\psi^0\left(-d\right)-\psi^0\left(d+\frac{n}{2}\right)\right)\qquad \textnormal{as } s\to 0, \label{eq: d not in N}
\end{align}
where $\Psi^0$ is the digamma function. $\qed$

In order to be able to use quasi-interpolation, the Fourier transform requires a singularity of even order at the origin. This is only achievable if the dimension $n$ is even and $d$ is an positive integer. So this seems to complement the case of the generalized multiquadric \cite{ortmann}.

\begin{theorem}
Let the dimension $n$ be even and $d$ and positive integer. Then the asymptotic behaviour up to and including the first logarithmic term is given by

\begin{align}
\hat\phi_\RN{3}(s)\sim&\sum_{j=0}^{1+\left\lfloor\frac{n}{2d}\right\rfloor}C_j(c,n,d)s^{-n+2d(j-1)}+\sum_{j=0}^{m_0-1}\tilde C_j(c,n,d)s^{2j}\\
&+\hat C_1(c,n,d) s^{2m_0}\left(\log\left(\frac{cs}{2}\right) +\hat C_2(c,n,d)\right), \qquad \textnormal{as } s\to 0,
\end{align}
where $m_0= d-(\frac{n}{2}\mod d)$ and the constants are given by

\begin{align}
C_0(c,n,d)&=(-1)^{d+1} \pi ^{n/2} d! 2^{2 d+n} \Gamma \left(d+\frac{n}{2}\right)\,,\\
C_j(c,n,d)&=\frac{(-1)^{j+1} \pi ^{n/2} 2^{2 d(1-j)+n}c^{2 d j} \Gamma \left(-j d+d+\frac{n}{2}\right)}{j \Gamma (j d-d)}\,, \qquad& j>0\\
\tilde C_j(c,n,d)&=\frac{(-1)^j\pi ^{n/2+1} 2^{-2j}  c^{2 d+2 j+n}}{ j! d \left(-\frac{n+2j}{2 d}-1\right) \Gamma \left(j+\frac{n}{2}\right) \sin \left(\pi  \left(-\frac{n+2j}{2 d}-1\right)\right)}\,,\qquad& j\geq 0\\
\hat C_1(c,n,d)&=\frac{ (-1)^{-\frac{m_0}{d}-\frac{n}{2 d}+m_0-1} 2^{1-2 m_0} \pi ^{n/2} d c^{2 d+2 m_0+n}}{m_0! \left(d+m_0+\frac{n}{2}\right) \Gamma \left(m_0+\frac{n}{2}\right)}\,,\\
\hat C_2(c,n,d)&=\frac{1}{2} \left(-\frac{2}{2 d+2 m_0+n}-H_{m_0}-\psi ^{(0)}\left(m_0+\frac{n}{2}\right)+\gamma \right).
\end{align}
$H_{m_0}=\sum_{k=1}^{m_0} \frac{1}{k}$ is the harmonic number and $\gamma$ is the Euler-Mascheroni constant.

\end{theorem}

\subsection*{\textit{Proof}}
Double poles appear when $t=-1-\frac{n}{2d}-\frac{m}{d}$ evaluates to a negative integer for some $m\in \mathbb{N}_0$. The first time this is given for $m_0=d-(n/2\mod d)$ and the double pole is then at $t_0=-1-\left \lceil\frac{n}{2d}\right\rceil$. All simple poles at negative integers up to $t_0$ are summarized in the first sum using the constants $C_j$. The second sum containing the constants $\tilde C_j$ comes from the simple poles at $-1-\frac{n}{2d},-1-\frac{n}{2d}-\frac{1}{d}, \dots, t_0+\frac{1}{d}$. The last term is derived by evaluating the second order pole at $t_0$.

\begin{theorem}
Let the dimension $n$ be even and $d$ a positive integer. Then the asymptotic behaviour of $\hat\phi_\RN{4}$ up to and including the first logarithmic term is given by

\begin{align}
\frac{\hat\phi_\RN{4}(s)}{2^n\pi^{n/2}c^{2d}}\sim&  \sum_{j=0}^{t_0-1}C_j(c,n,d)s^{-n+2dj}+\sum_{j=0}^{m_0-1}\tilde C_j(c,n,d)s^{2j}\\
&+\hat C_1(c,n,d) s^{2m_0}\left(\log\left(\frac{cs}{2}\right) +\hat C_2(c,n,d)\right)\,, \qquad \textnormal{as } s\to 0_+,
\end{align}
where $m_0= d-(\frac{n}{2}\mod d)$ as above. The constants are given by

\begin{align}
C_0(c,n,d)&=d \Gamma \left(\frac{n}{2}\right)\,,\\
C_j(c,n,d)&=\frac{(-1)^{j+1} \Gamma \left(\frac{n}{2}-dj\right)}{j \Gamma (dj)}\left(\frac{c}{2}\right)^{2dj}\,,\qquad & j>0\\
\tilde C_j(c,n,d)&=\frac{(-1)^j\pi }{ j!  \left(n/2+j\right) \Gamma \left(\frac{n}{2}+j\right) \sin \left(\pi  \left(\frac{n/2+j}{d}\right)\right)}\left(\frac{c}{2}\right)^{2j}\,,\qquad& j\geq 0\\
\hat C_1(c,n,d)&=\frac{ (-1)^{t_0+m_0} 2d}{m_0! \Gamma \left(1+m_0+\frac{n}{2}\right)}\,,\\
\hat C_2(c,n,d)&=\frac{1}{2} \left(\frac{1}{n/2+m_0}-H_{m_0}-\psi ^{(0)}\left(m_0+\frac{n}{2}\right)+\gamma \right).
\end{align}

\end{theorem}

\subsection*{\textit{Proof}}
The proof is analog to proof of the asymptotic behaviour of $\hat\phi_\RN{3}$, where the first double pole appears for $t=-\frac{n}{2d}-\frac{m}{d}$, where $m \in\mathbb{N}$ is the smallest integer such that $t$ is a negative integer. We call the smallest integer $m_0$ and the first double pole $t_0$. Like for $\hat\phi_{\RN{3}}$, $m_0$ is given by $m_0= d-(\frac{n}{2}\mod d)$ and $t_0=-\left \lceil\frac{n}{2d}\right\rceil$.
Adding the asymptotic behaviours of $\hat\phi_\RN{3}$ and $\hat\phi_\RN{4}$  results in an asymptotic series of the same structure as  $\hat\phi_\RN{3}$  with the sum of the constants from both series.

\begin{theorem}\label{theo: decayProperty} 

There exist coefficients $\mu_{\alpha ,\beta}$ such that the quasi-Lagrange function $\Psi(x):=\sum_{\beta \in \mathbb{Z}^n}\mu_{\beta}\phi(x-\beta)$ have the decay property
\begin{align}
\left |\Psi\right | \leq C \left(1+||x||\right)^{-2n-2d-2 m_0},\qquad x\in\mathbb{R}^n,
\end{align}
for some constant $C$ and furthermore, the operator defined by ("quasi-interpolation")
\begin{align}
Qf(x):=\sum_{j\in\mathbb{Z}^n} f(jh)\Psi\left(\frac{x-jh}{h}\right),\qquad x\in \mathbb{R}^n
\end{align}
reproduces all polynomials of degree $\ell_0=n+2d-1$ and converges uniformly with order
\begin{align}
||Qf-f||_\infty =\begin{cases}\mathcal{O}(h^{n+2d}\log(1/h))\qquad &m_0=0\\
\mathcal{O}(h^{n+2d})\qquad &m_0>0,
\end{cases}
\end{align}
where $0<h<1$ is the fill distance.
\end{theorem}

\subsection*{\textit{Proof}}
The proof employs Theorem A from \cite{BUHMANN2015156} and demonstrates that all requisite conditions are satisfied. In order to facilitate comparison, the notation from \cite{BUHMANN2015156}  is also employed.
Let $\hat \phi(\xi)=\frac{F(\xi)+bG_0(\xi)\log(\xi)}{G(\xi)}$ the asymptotic expansion as $\xi \to 0$. Comparing this with our previously derived asymptotic form we get

\begin{itemize}
\item $G_0(\xi)$ is a homogeneous polynomial of degree $n_0=n+2d+2 m_0$,
\item $G(\xi)$ is a homogeneous polynomial of degree $n+2d$ and $G(\xi)\neq 0$ for all $\xi \in \mathbb{R}^n\backslash\{0\}$,
\item $F(\xi)$ is a polynomial of degree $n_0$  with $F(0)\neq 0$. Therefore the Taylor series of $F$ at $\xi=0$ for some $\tilde m>n_0$ is error free.
\item The last property to show is that the Fourier transform $\hat \phi$ is continuously differentiable and that 
\begin{align}
\max_{|\gamma|\leq m_0+n+1}\int_{||\xi|| \geq \frac{1}{4}}|D^\gamma\hat\phi(\xi)| \, \mathrm{d}\xi\leq C<\infty.
\end{align}
\end{itemize}
Every derivative $D^\gamma$ can expressed in $n$-dimensional hyper-sphere coordinates and acting on some radial function is some linear combination of derivatives on the radial part. Therefore it is sufficient to show that the expression is valid for partial derivatives with respect to the radius at given order.

\begin{align}
\frac{\partial^\gamma}{\partial s^\gamma}\hat\phi_\RN{3}(s)=\frac{\partial^\gamma}{\partial s^\gamma} \frac{2^{n+2d}\pi^{\frac{n}{2}}} {s^{n+2d}} \frac{1}{2\pi \mathrm{i}} \underbrace{\int_L\frac{\pi}{t\sin(\pi t)}\frac{\Gamma\left(n/2+dt+d\right)}{\Gamma\left(-dt-d\right)}\left(\frac{2}{cs}\right)^{2dt}\mathrm{d}t}_{=:\Lambda}
\end{align}
Further use the parametrisation of $L$ by $t=r+i\tilde t$, where $r>0$ is fixed but arbitrarily and consider only the integral part.
\begin{align}
 \frac{\partial^\gamma}{\partial s^\gamma}  \Lambda&=\int_{-\infty}^\infty \frac{\pi}{(r+i \tilde t)\sin(\pi (r+i\tilde t))}\frac{\Gamma\left(n/2+dr+di\tilde t+d\right)}{\Gamma\left(-dr-di\tilde t-d\right)}\frac{\partial^\gamma}{\partial s^\gamma} \left(\frac{2}{cs}\right)^{2dr+2di\tilde t}\mathrm{d}\tilde t\\
&\leq \left(\frac{2}{cs}\right)^{2dr+\gamma}\int_{-\infty}^\infty \left |		\frac{\pi}{(r+i \tilde t)\sin(\pi (r+i\tilde t))}\frac{\Gamma\left(n/2+dr+di\tilde t+d\right)}{\Gamma\left(-dr-di\tilde t-d\right)}\frac{\Gamma\left( 2dr+2di\tilde t+\gamma\right)}{\Gamma\left( 2dr+2di\tilde t\right)} 	\right|\mathrm{d}\tilde t.	
\end{align}
To show that the integral is finite for all $r>0$ we look the asymptotic behaviour of the integrand. Therefore we use the formulas given in \cite{abr}
\begin{align}
\left |\Gamma \left(x+iy\right)\right | \sim \sqrt{2\pi} |y|^{x-1/2}\textnormal{e}^{-\pi|y|/2}\qquad \textnormal {as } y\to \pm\infty,
\end{align}
where $x$ is a bounded real value,

\begin{align}
\left|\frac{\Gamma\left( 2dr+2di\tilde t+\gamma\right)}{\Gamma\left( 2dr+2di\tilde t\right)} \right| &\sim |2d\tilde t|^\gamma, \qquad &\textnormal {as } \tilde t\to \pm\infty,\\
\left|\frac{\Gamma\left(n/2+dr+di\tilde t+d\right)}{\Gamma\left(-dr-di\tilde t-d\right)}\right| &\sim |d\tilde t|^{\frac{n}{2}+2dr+2d},\qquad &\textnormal {as } \tilde t\to \pm\infty,\\
r+i\tilde t&\sim i\tilde t  ,&\textnormal {as } \tilde t\to \pm\infty,\\
|\sin\left( \pi r+i\pi \tilde t\right)| &\sim\frac{\textnormal{e} ^{\pi |\tilde t|}}{2} , &\textnormal {as } \tilde t\to \pm\infty.
\end{align}
In total it is straightforward to see that the integrated function decays faster than any polynomial and so for all $r>0$ and for all $\gamma$ we can evaluate the integral to a constant.
Since $r>0$ can be arbitrarily large, the Fourier transform $\hat \phi$ and all its derivatives decays faster than every polynomial and so
\begin{align}
\max_{|\gamma|\leq m_0+n+1}\int_{||\xi||\geq\frac{1}{4}}|D^\gamma\hat\phi(\xi)|\mathrm{d}\xi \leq C<\infty
\end{align}
 is finite. Analog statements are true for $\hat\phi_\RN{4}$. In the classical case $d=1$, the Fourier transform is given by some power of $s$ times the modified Besselfunktion. For large arguments the modified Bessel function decays exponential and so all derivatives do. So this result is in perfekt agreement.

\section{Shift-invariant subspaces of Sobolev spaces}
Following the notation in \cite{article} the Sobolev space of smoothness is defined as
\begin{align}
W_2^s(\mathbb{R}^n):=\{f\in \mathcal{S}'(\mathbb{R}^n):\int_{\mathbb{R}^n}|\mathcal{F}f(x)|^2 (1+||x||)^s \mathrm{d}x <\infty\},
\end{align}
where $s\in \mathbb{R}$. The used norm is given by
\begin{align}
||f||_{W_2^s(\mathbb{R}^n)}:= \left( \int_{\mathbb{R}^n}|\mathcal{F}f(x)|^2 (1+||x||)^s \mathrm{d}x\right)^{\frac{1}{2}}.
\end{align}
Let $S$ be a \textit{principal shift-invariant space} then we define $S^h:=\{f(\cdot/h)| f\in S\}$, $ h>0$ and the sequence $\mathcal{S}:=\{ S^h: h>0\}$ is called a ladder. Furthermore the ladder provides approximation oder $k$ in $W_2^s(\mathbb{R}^n)$ if there exists a constant $C$ such that for every $f \in W_2^s(\mathbb{R}^n)$ and for every $h>0$
\begin{align}
\inf_{g\in S^h}||f-g||_{W_2^s(\mathbb{R}^n)}\leq C h^{k-s} ||f||_{W_2^s(\mathbb{R}^n)}.
\end{align}
The following theorem is shown in \cite{article}
\begin{theorem}\label{sec: Subspaces}
Let $k\geq0, s<k, \Psi \in W_2^s(\mathbb{R}^n)$. Suppose that, for some $\eta_1,\eta_2 >0$ and for some ball $\Omega$ centered at the origin,
\begin{align}
\eta_1\leq |\hat\Psi|\leq \eta_2 \qquad a.e.\, on \, \Omega.
\end{align}
Then $S_\Psi (W_2^s(\mathbb{R}^n))$ provides approximation order k in $W_2^s(\mathbb{R}^n)$ if and only if

\begin{itemize}
\item[a)] For all $\beta \in 2\pi\mathbb{Z}^n\backslash \{0\}$ $$c_\beta :=\lim_{||\xi|| \to 0}||\xi||^{-k} |\hat\Psi(\xi+\beta)|<\infty$$  
\item[b)] and $$\sum_{\beta \in 2\pi\mathbb{Z}^n\backslash\{0\}} c_\beta^2\cdot ||\beta||^{2s}< \infty.$$
\end{itemize}
\end{theorem}

Now we check the approximation order of our constructed quasi-Lagrange function $\Psi$ using Theorem~\ref{sec: Subspaces}.
Since  $\hat\Psi$ is continuous and $\hat\Psi(0)=1$,  there exists two constants $\eta_1,\eta_2>0$ such that $\hat\Psi$ is bounded on some open ball centered at the origin.
To check Part $a)$ we rewrite the Fourier transform of the quasi-Lagrange function as
\begin{align}
\hat\Psi(x):=\sum_{\alpha\in\mathbb{Z}^n}\mu_\alpha \textnormal{e}^{-i\alpha\cdot x} \hat\phi(x)=P(x) \hat\phi(x).
\end{align}
Since $\alpha$ is an integer-touple, the function $P(x)$ is $2\pi$-periodic in every dimension and
\begin{align}
\lim_{||\xi|| \to 0}||\xi||^{-k}|\hat\Psi(\xi+\beta)|= \lim_{\xi\to 0}\frac{|P(\xi)|}{||\xi||^k}|\hat\phi(\beta)|.
\end{align}
The Taylor expansion of $P(\xi)$ at zero matches the order of the singularity of the Fourier transform of our RBF. It follows
\begin{align}
c_\beta=\begin{cases}
0 \qquad &k<n+2d\\
C |\hat\phi(\beta)| \qquad &k=n+2d\\
\infty \qquad &k>n+2d,
\end{cases}
\end{align}
for some constant $C$. The only interesting case is for $k=n+2d$, the others are trivial. The only singularity of $\hat\phi$ is at the origin and since $\beta\neq 0$ the expression evaluates to an finite number. Furthermore the decay rate of $c_\beta$ is given by the decay rate of $\hat\phi$. In the proof of Theorem~\ref{theo: decayProperty}, we showed that $\hat\phi$ decays faster than every polynomial and so it follows that the sum in Part $b)$ is finite for all $s<k$. The approximation order of the constructed quasi-Lagrange function is only affected by the singularity at the origin and provides approximation order $n+2d$.

Comparing this result in our case with the previous one by Strang and Fix, it is evident that quasi-interpolation achieves the highest attainable approximation order. In this case quasi-interpolation provides an constructive way of getting the best approximation in the spanned principal shift-invariant subspace.


\section{New generalized multiquadric- and inverse multiquadric function}
This section introduces a further generalization of multi- and inversmultiquadrics within the framework of quasi-interpolation. The RBF under analysis is 
\begin{align}
\phi(x)=\left ( c^\lambda+||x||^\lambda \right )^\beta, \qquad c>0, x\in\mathbb{R}^n, \lambda\in\mathbb{R},\beta\in\mathbb{R}\backslash \mathbb{N}.
\end{align}
Classical cases are the Hardy-multiquadric for $\lambda=2, \beta =\frac{1}{2}$, the inverse multiqudrics for $\lambda=2, \beta =-\frac{1}{2}$,  the generalized multiquadric $\lambda =2$, $\beta$ arbitrarily and Ortmann and Buhmann \cite{ortmann} proposed $\lambda =2d, \beta =\frac{1}{2}$, where $d\in \mathbb{N}$. As in our previous analysis, we commence with the calculation of the Fourier transform. Subsequently, we will examine the asymptotic behaviour in different scenarios.

\begin{theorem}\label{FT2}
The generalized $n$-dimensional Fourier transform of
 \begin{align}
\phi(x)=\left ( c^\lambda+||x||^\lambda \right )^\beta, \qquad  c>0,\,  \,x \in \mathbb{R}^n, r=||x||, \lambda\in\mathbb{R},\beta\in\mathbb{R}\backslash \mathbb{N}\\
\end{align}
is given by the Mellin-Barnes type integral
\begin{align}
\hat\phi (s)=\frac{2^{n}\pi^{\frac{n}{2}}c^{\lambda\beta}} {s^{n}\Gamma\left(-\beta\right)} \frac{1}{2\pi \mathrm{i}} \int_L\frac{\Gamma \left(-t\right)\Gamma \left(-\beta+t\right)\Gamma\left( \frac{n+\lambda t}{2}\right)}{\Gamma\left(-\frac{\lambda t}{2}\right)}\left(\frac{2}{cs}\right)^{\lambda t}\mathrm{d}t \label{eq: FT IMQ},
\end{align}
where $L$ is the path coming from $-\mathrm{i}\infty $ to $\mathrm{i}\infty$, while the poles of Gamma functions in the numerator are separed.
\end{theorem}
\subsection*{\textit{Proof}}
The proof follows the same steps as shown above respectivly in \cite{ortmann} with $\lambda=2d$. $\beta=\frac{1}{2}$ and
\begin{align}
\left ( c^\lambda+r^\lambda \right )^\beta=c^{\lambda\beta}G_{1,1}^{\,1,1}\!\left(\left.{\begin{matrix}1+\beta\\ 0\end{matrix}}\;\right|\,\left(\frac{r}{c}\right )^{\lambda}\right).\qed
\end{align}

\subsection{Asymptotic behaviour}
In order to analyse the asymptotic behaviour of the Fourier transform for $s \to 0_+$, it is necessary to consider a number of different cases. In all cases, a simple pole leads to the presence of powers of $s$, whereas higher-order poles leads to logarithmic singularities. 
In the context of high order quasi-interpolation, it is preferable to have as many powers of $s$ present as possible before the first logarithmic singularity appears. In addition, we will perform the calculations for a further intriguing case.
\subsubsection{Generalized multiquadrics $\lambda,\beta >0$}
\begin{figure}
\centering
\input{VisualisierungPolstellen1.tex}
\caption{Visualization of the poles from equation~(\ref{eq: FT IMQ}) for $\lambda>0, \beta>0$. The singularities of the gamma functions in the nominator are depicted as dots, whereas the singularities of the gamma function in the denominator are represented by circles (or holes), which have the potential to cancel poles. For the purposes of enhanced visualisation, the poles have been represented as if they had an imaginary component, although this is not the case. Every pole lies on the real axis. Additionally, the path of integration is depicted and indicates which poles must be included. The example is provided for $|\lambda|=2, |\beta|=\frac{3}{2}, n=1$.  }
\label{fig: Visualisierung1}
\end{figure}

The poles are visualized in figure \ref{fig: Visualisierung1}. The order of the poles are given in descending order from $\beta$ to $-\infty$. For quasi-interpolation the most useful case is $\beta \notin\mathbb{N}$ with $\lambda \beta$ is not even. Then the asymptotic behaviour is given by 
 \begin{align}
\hat\phi(s)\sim \frac{2^{n+\lambda\beta}\pi^{\frac{n}{2}}}{s^{n+\lambda\beta}}\frac{\Gamma\left(\frac{n+\lambda\beta}{2}\right)}{\Gamma\left(-\frac{\lambda\beta}{2}\right)},\qquad \textnormal{as } s\to 0_+.
\end{align}
The special case $\lambda =2d$ and $\beta=\frac{1}{2}$, $d \in \mathbb{N}$ with all calculations of the asymptotic behaviour can be found in \cite{ortmann}. They used only odd dimensions to get a trignonometric polynomial with finitely many coeffients. To get as many poles between $\beta$ and the first double pole they increased $\lambda$. But in figure~\ref{fig: Visualisierung1} it is easy to see that the same thing can be done increasing $\beta$. Results in the same fashion could also be achieved using for example $\lambda = 3$ and $\beta =\frac{1}{3}$. If the leading order is not of even power it is also possible to use the function in quasi-interpolation with a trigonometric polynomial using infinitely many coeffients \cite{Buhmann2023}.

\subsubsection{Classic generalized inverse multiquadrics $\lambda >0 ,\beta <0$}
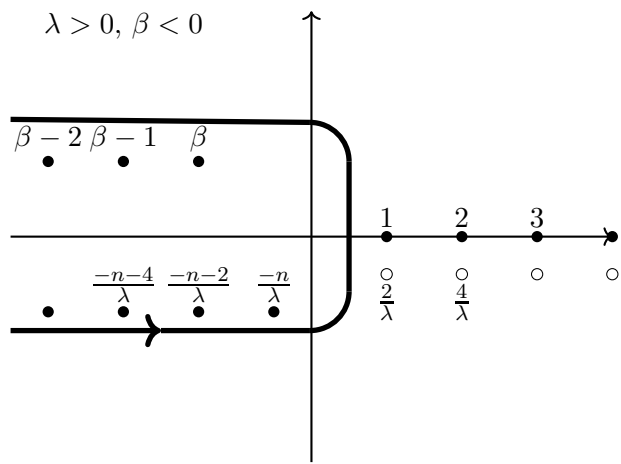
\begin{figure}
\centering
\input{VisualisierungPolstellen2.tex}
\caption{Visualization of the poles from equation~(\ref{eq: FT IMQ}) for $\lambda>0, \beta<0$. The singularities of the gamma functions in the nominator are depicted as dots, while the singularities of the gamma function in the denominator are represented by circles (or holes), which can cancel poles. For purposes of enhanced visualisation, the poles have been represented as if they had an imaginary component, although this is not the case. Every pole lies on the real axis. Additionally, the path of integration is depicted and indicates which poles must be included. The example is provided for $|\lambda|=2, |\beta|=\frac{3}{2}, n=1$.  }
\label{fig: Visualisierung2}
\end{figure}

The poles for this case are illustrated in figure~\ref{fig: Visualisierung2}. Let $\tilde\beta =-\beta >0$ and $\tilde\beta > \frac{n}{\lambda}$ then the first simple pole is at $t=\frac{n}{\lambda}$. Evaluating this pole leads to leading asymptotic behaviour of a constant. Consequently, quasi-interpolation cannot be employed in this case. If $\beta = \frac{n}{\lambda}$ like for $n=1$ the classical inverse multiquadric ($\beta=-\frac{1}{2}, \lambda=2$) the leading asymptotic behaviour is given by a logarithmic term. Buhmann showed that with an infinite expansion of a trigonometric polynomial this can also be used in quasi-interpolation \cite{Buhmann1993}. The final case $\tilde \beta < \frac{n}{\lambda}$ is of particular interest. The first pole is located at $t = \beta$, and evaluating the residue leads to the following result: 
\begin{align}
\hat\phi(s)\sim \frac{2^{n-\lambda\tilde\beta}\pi^{\frac{n}{2}}}{s^{n-\lambda\tilde\beta}}\frac{\Gamma\left(\frac{n-\lambda\tilde\beta}{2}\right)}{\Gamma\left(\frac{\lambda\tilde \beta}{2}\right)},\qquad \textnormal{as } s\to 0_+.
\end{align}
If the condition $n>\lambda\beta$ is met, a singularity occurs at the origin. In this case, the function can be used for quasi-interpolation. Otherwise, the function cannot be used for quasi-interpolation.

\subsubsection{Singular RBF $\lambda <0 ,\beta >0$}
\begin{figure}
\centering
\input{VisualisierungPolstellen3.tex}
\caption{Visualization of the poles from equation~(\ref{eq: FT IMQ}) for $\lambda<0, \beta>0$. The singularities of the gamma functions in the nominator are depicted as dots, while the singularities of the gamma function in the denominator are represented by circles (or holes), which can cancel poles. For purposes of enhanced visualisation, the poles have been represented as if they had an imaginary component, although this is not the case. Every pole lies on the real axis. Additionally, the path of integration is depicted and indicates which poles must be included. The example is provided for $|\lambda|=2, |\beta|=\frac{3}{2}, n=1$.  }
\label{fig: Visualisierung3}
\end{figure}
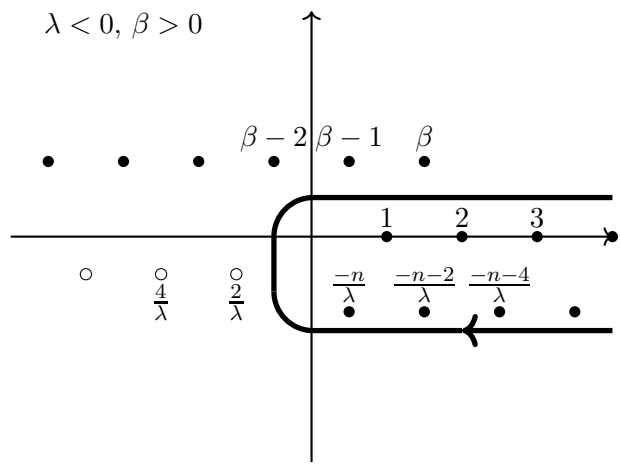
The poles for this case are shown in figure \ref{fig: Visualisierung3}. For $\lambda <0$ the Mellin-Barnes-integral changes the integral contour. Let $\tilde\lambda =-\lambda >0$. If $\beta = \frac{n}{\tilde\lambda}$, then the leading order is given by the evaluation of a double pole, which leads to a logarithmic singularity of the Fourier transform at the origin. Further if $ \frac{n}{\tilde \lambda}< 1$ then the first pole is at $t = \frac{n}{\tilde \lambda}$ and evaluates to a constant. Consequently, this case is not suitable for quasi-interpolation. Finally, if  $\frac{n}{\tilde\lambda}>1$, the first pole is located at $t=1$, and its evaluation leads to
\begin{align}
\hat\phi(s)\sim \frac{2^{n-\tilde\lambda}\pi^{\frac{n}{2}}c^{\tilde \lambda(1-\beta)}\beta} {s^{n-\tilde\lambda}} \frac{\Gamma\left( \frac{n-\tilde\lambda }{2}\right)}{\Gamma\left(\frac{\tilde\lambda}{2}\right)} ,\qquad \textnormal{as } s\to 0_+.
\end{align}
So again there is a singularity which can be used for quasi-interpolation. Remember that the RBF itself is singular at the origin. That means that this can not be used directly. Notice that this is the first time the leading order depends on the constant $c$. To use this RBF one have to take a linear combination using different constants $c$. For example $\tilde\lambda=2, \beta=\frac{1}{2}$
 $$\phi(r)=\left(\frac{1}{c_1^{2}}+\frac{1}{r^{2}}\right)^{\frac{1}{2}}-\left(\frac{1}{c_2^{2}}+\frac{1}{r^{2}}\right)^{\frac{1}{2}}\sim \frac{1}{2} r \left(\frac{1}{c_1^2}-\frac{1}{c_2^2}\right)\qquad \textnormal{as } r\to 0_+ .$$ 
To cancel out the singularites of the RBF for $\beta>1$ more terms are needed with specific conditions to the constants $c$. Nevertheless, there remains one additional degree of freedom, which allows the linear combination to avoid cancelling out the highest order of the singularity of the Fourier transform. For $n\geq 3$ it is possible to choose $\tilde\lambda$ such that a finite linear combination for a trigonometric polynomial can be choosen. Consequently, this RBF has the potential to overcome the curse of dimensionality. This is a remarkable result. However, the analysis is not part of this paper and will be part of further research.

\subsubsection{Inverse multiquadric- type RBF $\lambda <0 ,\beta <0$}
\begin{figure}
\centering
\input{VisualisierungPolstellen4.tex}
\caption{Visualization of the poles from equation~(\ref{eq: FT IMQ}) for $\lambda<0, \beta<0$. The singularities of the gamma functions in the nominator are depicted as dots, while the singularities of the gamma function in the denominator are represented by circles (or holes), which can cancel poles. For purposes of enhanced visualisation, the poles have been represented as if they had an imaginary component, although this is not the case. Every pole lies on the real axis. Additionally, the path of integration is depicted and indicates which poles must be included. The example is provided for $|\lambda|=2, |\beta|=\frac{3}{2}, n=1$.  }
\label{fig: Visualisierung4}
\end{figure}
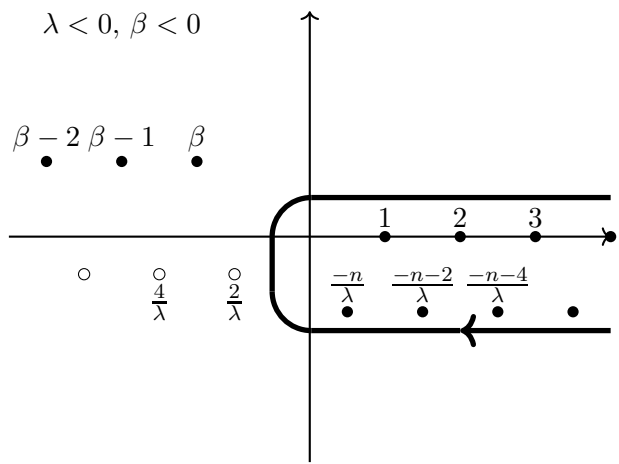
The poles of this case are shown in figure \ref{fig: Visualisierung4}. To the best of our knowledge, this case has not been considered previously, either as a general nor as a special case. Again let $\tilde \lambda=-\lambda >0$ and $\tilde \beta =-\beta >0$. As before $\frac{n}{\tilde \lambda}$ should be greater than one because then the Fourier transform has a singularity at the origin. If so the asymptotic behaviour is given by 
\begin{align}
\hat\phi(s)\sim- \frac{2^{n-\tilde\lambda}\pi^{\frac{n}{2}}c^{\tilde \lambda(1+\tilde \beta)}\tilde \beta} {s^{n-\tilde\lambda}} \frac{\Gamma\left( \frac{n-\tilde\lambda }{2}\right)}{\Gamma\left(\frac{\tilde\lambda}{2}\right)} ,\qquad \textnormal{as } s\to 0_+.
\end{align}
The only difference in the Fourier transform of the case of the singular RBF is that $\beta$ ( which does only affect the factor) is negative and that the RBF itself is not singular. Sono linear combination of RBF is required for quasi-interpolation. Again for $n\geq 3$ the singularity at the origin is affected by $\lambda$ and thus the curse of dimensionality is overcome. Consequently, we undertake a more detailed analysis of this case, specifically focusing on determining the highest approximation order of the quasi-interpolant.

Evaluating the simple poles for some values of $n, \beta, \lambda$ such that no double pole appears leads to
\begin{align}
\hat\phi(s)&= \sum_{k=1}^\infty (-1)^k \frac{ 2^{n-\tilde\lambda k}\pi^{\frac{n}{2}}c^{\tilde \lambda(k-\tilde \beta)}} {k! s^{n-\tilde\lambda k}} \frac{\Gamma\left(\tilde\beta+k\right)}{\Gamma\left(\tilde\beta \right)}
\frac{\Gamma\left( \frac{n-\tilde\lambda k }{2}\right)}{\Gamma\left(\frac{\tilde\lambda k}{2}\right)}\\
&+\sum_{k=0}^\infty (-1)^k \frac{ \pi^{\frac{n}{2}}c^{\tilde \lambda\tilde \beta+n+2k}} {2^{2 k-1}\tilde\lambda k!}
\frac{\Gamma\left( -\frac{n+2k }{\tilde\lambda}\right)\Gamma\left( \tilde\beta +\frac{n+2k }{\tilde\lambda}\right)}{\Gamma\left( \tilde\beta\right) \Gamma\left(\frac{n}{2}+k\right)}s^{2k}.		\label{eq: Asymptotic Expansion}
\end{align}
It should be noted that in the first sum of equation~(\ref{eq: Asymptotic Expansion}) odd powers of $s$ may appear. 
Should one wish to perform classical quasi-interpolation, a double pole will emerge, resulting in the appearance of a logarithmic term or an odd power of s. Consequently, the approximation order will be constrained.

\begin{theorem}
Let $\beta <0$ and the dimension $n\geq 3$ then there exists $\lambda <0$, $\beta <0$ and finitely many coeffients $\mu_k$  such that 
$$\Psi(x) =\sum_{k\in\mathbb{Z}^n} \mu_k \phi(||x-k||),\qquad x\in\mathbb{R}^n$$ 
reproduces all polynomials of degree $\frac{n-1}{2}$  if $n$ is odd and of degree $n-3$ if $n$ is even and the uniform approximation error can be estimated by
\begin{align}
\Vert Q_hf-f\Vert_{\infty}=\begin{dcases}
\mathcal{O}\left(h^{n-2}\log (1/h) \right)\qquad &n\textnormal{ even}\\
\mathcal{O}\left(h^{\frac{n-3}{2}} \right)\qquad &n \textnormal{ odd.}\\
\end{dcases}
\end{align}
\end{theorem}

\subsection*{\textit{Proof}}

As in the first part of this paper we use Theorem A from \cite{BUHMANN2015156}.  In order to achieve this, we undertake a comparison of the asymptotic behaviour of our RBF with that of the function $\hat\phi (\xi) = \frac{F(\xi)+bG_0(\xi)\log(|\xi|)}{G(\xi)}$ as $\xi \to 0_+$. This enables us to demonstrate that the criteria set out in the theorem are satisfied.
\subsubsection*{$n$ is even}
Let $n\geq 3$ be even and $\tilde \lambda=-\lambda=2$. The first double pole appears for $t=\frac{n}{2}$. So we get
\begin{align}
\hat\phi(s)\sim\sum_{k=1}^{\frac{n}{2}}C_ks^{-n+2k}+\tilde C \log(s),\qquad \textnormal{as } s\to 0_+
\end{align}
for some  constants $C_k$ and $\tilde C$. Since all powers of $s$ are even there exists a finite trigonometric polynomial, which cancels them out. The points 2. and 3. of the Strang and Fix conditions are satisfied for $m=n-3$. Furthermore $F(\xi)$ is a polynomial so it is smooth and the Taylor expansion is exact. The asymptotic behaviour of $\hat\phi (s) $ for $s\to \infty$ is limited by the pole at $\beta$ if $\tilde \beta \notin \mathbb{N}$. The decay rate can be determined by employing the specified path $L=\beta+\varepsilon +\mathrm{i} T, T\in (-\infty, \infty)$ for arbitrarily small $\varepsilon>0$. Hence,
\begin{align}
\hat\phi(s)\sim C s^{-n-\lambda\beta},\qquad \textnormal{as }s\to \infty,
\end{align}
with some constant $C$. Higher derivatives decay even faster as $s\to \infty$.
For $\tilde\beta \in\mathbb{N}$ there is no pole on the left half plane that limitates  the contour integral. It follows that $\hat\phi(s)$ then decays faster than every polynomial and the same applies for all derivatives. So for both cases  $\int_{||\xi||>\frac{1}{4}}  |D^\gamma \hat\phi(\xi)| \mathrm{d}\xi <\infty \quad\forall \gamma\in\mathbb{Z}^n $ is finite and using Theorem~A from \cite{BUHMANN2015156} we get
\begin{align}
|\Psi (x)|\leq C(1+||x||)^{-2n+2},\qquad x\in\mathbb{R}^n.
\end{align}
Using the Theorem of Strang and Fix we get the desired result.

\subsubsection*{$n$ is odd}
Let $n\geq 3$ be odd and $\tilde \lambda=-\lambda=\frac{n-1}{2}$.
If $\tilde \lambda$ is even there is no double pole and the asymptotic behaviour is fully given in equation~(\ref{eq: Asymptotic Expansion}). If $\tilde\lambda $ is odd, we end up with
\begin{align}
\hat\phi(s)\sim C_1 s^{-\frac{n+1}{2}} +C_2s^{-1}+C_3 s^{\frac{n-3}{2}}\log(s)+\sum_{k=0}^{\frac{n-3}{4}}\tilde C_ks^{2k},\qquad \textnormal{as } s\to 0_+
\end{align}
for some constants $C_1, C_2, C_3$ and $\tilde C_k$. In both cases the second summand ($k=2$) in the first sum of equation~(\ref{eq: Asymptotic Expansion}) results in an odd power of $s$. Note that the second sum only introduces even powers of $s$ which are not necessary for the argumentation because a finite trigonometric polynomial can handle them. So the points 2. and 3. of the Strang and Fix conditions are satisfied for $m=\frac{n-3}{2}$. The decay rate of $\hat\phi(s)$ follows the same argumentation as in the case $n$ even. It follows that 
$$|\Psi(x)| \leq C(1+||x||)^{-2n+1},\qquad x\in\mathbb{R}^n.$$ Using this in the Theorem of Strang and Fix we get the desired result.

%
%

\bibliographystyle{abbrv}
\bibliography{Quellen}

\end{document}

%% file: VisualisierungPolstellen1.tex
%
%
\begin{tikzpicture}

\tikzmath{%
\s=1;
\b=1.5;
\d=1;
\l=2;
\r=0.52;
}

\draw[thick,->] (-4,0) -- (4,0);
\draw[thick,->] (0,-3) -- (0,3);

\tkzDefPoint(-2.5,2.5){node1}
\tkzLabelPoint[above](node1){$\lambda>0,\,\beta>0$}

\tkzDefPoint(\b,1){b1}
\tkzDefPoint(\b-\s,1){b2}
\tkzDefPoint(\b-2*\s,1){b3}
\tkzLabelPoint[above](b1){$\beta$}
\tkzLabelPoint[above](b2){$\beta-1$}
\tkzLabelPoint[above](b3){$\beta-2$}

\foreach \n in {0,...,5}
	\node at (\b-\n*\s,1)[circle,fill,inner sep=1.5pt]{};

\tkzDefPoint(\s,0){t1}
\tkzDefPoint(2*\s,0){t2}
\tkzDefPoint(3*\s,0){t3}
\tkzLabelPoint[above](t1){$1$}
\tkzLabelPoint[above](t2){$2$}
\tkzLabelPoint[above](t3){$3$}

\foreach \n in {1,...,4}
	\node at (\n*\s,0)[circle,fill,inner sep=1.5pt]{};

\tkzDefPoint(-\d/\l*\s,-1){n1}
\tkzDefPoint(-\d/\l*\s-2/\l*\s,-1){n2}
\tkzDefPoint(-\d/\l*\s-4/\l*\s,-1){n3}
\tkzLabelPoint[above](n1){$\frac{-n}{\lambda}$}
\tkzLabelPoint[above](n2){$\frac{-n-2}{\lambda}$}
\tkzLabelPoint[above](n3){$\frac{-n-4}{\lambda}$}

\foreach \n in {1,...,4}
	\node at (\d/\l*\s-2/\l*\s*\n,-1)[circle,fill,inner sep=1.5pt]{};

\tkzDefPoint(2/\l*\s*1,-0.5){l2}
\tkzDefPoint(2/\l*\s*2,-0.5){l3}
\tkzLabelPoint[below](l2){$\frac{2}{\lambda}$}
\tkzLabelPoint[below](l3){$\frac{4}{\lambda}$}

\foreach \n in {1,...,4}
	\node at (2/\l*\s*\n,-0.5)[circle,draw,inner sep=1.5pt]{};

\coordinate (A) at (-4,-1.25);
\coordinate (B) at ($(A)+ (2,0)$) ;
\coordinate (C) at ($(A)+ (4.5,0)-(\r,0)$);
\coordinate (D) at ($(C)+(\r,\r)$);
\coordinate (E) at (0.5,0);
\coordinate (F) at  ($(E)+ (\r,\r)$) ;
\coordinate (G) at  ($(F)+ (\b-0.5-\r,0)$) ;
\coordinate (H) at (\b,3*\r);
\coordinate (J) at (-4,3*\r);

\coordinate (K) at (A)++(1,1);

\draw [black, line width=2pt, ->] plot [smooth] coordinates {(A) (B) };
\draw [black, line width=2pt] plot [smooth] coordinates {(B) (C) };
\draw[black, line width=2pt] (C) arc (-90:0:\r) ;
\draw [black, line width=2pt] plot [smooth] coordinates {(D)  (E)  };
\draw[black, line width=2pt] (E) arc (180:90:\r) ;
\draw [black, line width=2pt] plot [smooth] coordinates {(F) (G)  };
\draw[black, line width=2pt] (G) arc (-90:90:\r) ;
\draw [black, line width=2pt] plot [smooth] coordinates {(H) (J)  };

\end{tikzpicture}%

%% file: VisualisierungPolstellen2.tex
%
%
\begin{tikzpicture}

\tikzmath{%
\s=1;
\b=-1.5;
\d=1;
\l=2;
\r=0.52;
}

\draw[thick,->] (-4,0) -- (4,0);
\draw[thick,->] (0,-3) -- (0,3);

\tkzDefPoint(-2.5,2.5){node1}
\tkzLabelPoint[above](node1){$\lambda>0,\,\beta<0$}

\tkzDefPoint(\b,1){b1}
\tkzDefPoint(\b-\s,1){b2}
\tkzDefPoint(\b-2*\s,1){b3}
\tkzLabelPoint[above](b1){$\beta$}
\tkzLabelPoint[above](b2){$\beta-1$}
\tkzLabelPoint[above](b3){$\beta-2$}

\foreach \n in {0,...,2}
	\node at (\b-\n*\s,1)[circle,fill,inner sep=1.5pt]{};

\tkzDefPoint(\s,0){t1}
\tkzDefPoint(2*\s,0){t2}
\tkzDefPoint(3*\s,0){t3}
\tkzLabelPoint[above](t1){$1$}
\tkzLabelPoint[above](t2){$2$}
\tkzLabelPoint[above](t3){$3$}

\foreach \n in {1,...,4}
	\node at (\n*\s,0)[circle,fill,inner sep=1.5pt]{};

\tkzDefPoint(-\d/\l*\s,-1){n1}
\tkzDefPoint(-\d/\l*\s-2/\l*\s,-1){n2}
\tkzDefPoint(-\d/\l*\s-4/\l*\s,-1){n3}
\tkzLabelPoint[above](n1){$\frac{-n}{\lambda}$}
\tkzLabelPoint[above](n2){$\frac{-n-2}{\lambda}$}
\tkzLabelPoint[above](n3){$\frac{-n-4}{\lambda}$}

\foreach \n in {1,...,4}
	\node at (\d/\l*\s-2/\l*\s*\n,-1)[circle,fill,inner sep=1.5pt]{};

\tkzDefPoint(2/\l*\s*1,-0.5){l2}
\tkzDefPoint(2/\l*\s*2,-0.5){l3}
\tkzLabelPoint[below](l2){$\frac{2}{\lambda}$}
\tkzLabelPoint[below](l3){$\frac{4}{\lambda}$}

\foreach \n in {1,...,4}
	\node at (2/\l*\s*\n,-0.5)[circle,draw,inner sep=1.5pt]{};

\coordinate (A) at (-4,-1.25);
\coordinate (B) at ($(A)+ (2,0)$) ;
\coordinate (C) at ($(A)+ (4.5,0)-(\r,0)$);
\coordinate (D) at ($(C)+(\r,\r)$);
\coordinate (E) at (0.5,1);
\coordinate (F) at  ($(E)+ (-\r,\r)$) ;
\coordinate (G) at  ($(F)+ (\b-0.5-\r,0)$) ;
\coordinate (H) at (\b,3*\r);
\coordinate (J) at (-4,3*\r);

\coordinate (K) at (A)++(1,1);

\draw [black, line width=2pt, ->] plot [smooth] coordinates {(A) (B) };
\draw [black, line width=2pt] plot [smooth] coordinates {(B) (C) };
\draw[black, line width=2pt] (C) arc (-90:0:\r) ;
\draw [black, line width=2pt] plot [smooth] coordinates {(D)  (E)  };

\draw[black, line width=2pt] (E) arc (0:90:\r) ;
\draw [black, line width=2pt] plot [smooth] coordinates {(F) (J)  };

\end{tikzpicture}%

%% file: VisualisierungPolstellen3.tex
%
%
\begin{tikzpicture}

\tikzmath{%
\s=1;
\b=1.5;
\d=1;
\l=2;
\r=0.52;
}

\draw[thick,->] (-4,0) -- (4,0);
\draw[thick,->] (0,-3) -- (0,3);

\tkzDefPoint(-2.5,2.5){node1}
\tkzLabelPoint[above](node1){$\lambda<0,\,\beta>0$}

\tkzDefPoint(\b,1){b1}
\tkzDefPoint(\b-\s,1){b2}
\tkzDefPoint(\b-2*\s,1){b3}
\tkzLabelPoint[above](b1){$\beta$}
\tkzLabelPoint[above](b2){$\beta-1$}
\tkzLabelPoint[above](b3){$\beta-2$}

\foreach \n in {0,...,5}
	\node at (\b-\n*\s,1)[circle,fill,inner sep=1.5pt]{};

\tkzDefPoint(\s,0){t1}
\tkzDefPoint(2*\s,0){t2}
\tkzDefPoint(3*\s,0){t3}
\tkzLabelPoint[above](t1){$1$}
\tkzLabelPoint[above](t2){$2$}
\tkzLabelPoint[above](t3){$3$}

\foreach \n in {1,...,4}
	\node at (\n*\s,0)[circle,fill,inner sep=1.5pt]{};

\tkzDefPoint(\d/\l*\s,-1){n1}
\tkzDefPoint(\d/\l*\s+2/\l*\s,-1){n2}
\tkzDefPoint(\d/\l*\s+4/\l*\s,-1){n3}
\tkzLabelPoint[above](n1){$\frac{-n}{\lambda}$}
\tkzLabelPoint[above](n2){$\frac{-n-2}{\lambda}$}
\tkzLabelPoint[above](n3){$\frac{-n-4}{\lambda}$}

\foreach \n in {1,...,4}
	\node at (-\d/\l*\s+2/\l*\s*\n,-1)[circle,fill,inner sep=1.5pt]{};

\tkzDefPoint(-2/\l*\s*1,-0.5){l2}
\tkzDefPoint(-2/\l*\s*2,-0.5){l3}
\tkzLabelPoint[below](l2){$\frac{2}{\lambda}$}
\tkzLabelPoint[below](l3){$\frac{4}{\lambda}$}

\foreach \n in {1,...,3}
	\node at (-2/\l*\s*\n,-0.5)[circle,draw,inner sep=1.5pt]{};

\coordinate (A) at (4,-1.25);
\coordinate (B) at ($(A)+ (-2,0)$) ;
\coordinate (C) at ($(A)+ (-4.5,0)+(\r,0)$);
\coordinate (D) at ($(C)+(-\r,\r)$);
\coordinate (E) at (-0.5,0);
\coordinate (F) at  ($(E)+ (\r,\r)$) ;
\coordinate (G) at  ($(F)+ (\b-0.5-\r,0)$) ;
\coordinate (H) at (\b,3*\r);
\coordinate (J) at (4,\r);

\coordinate (K) at (A)++(1,1);

\draw [black, line width=2pt, ->] plot [smooth] coordinates {(A) (B) };
\draw [black, line width=2pt] plot [smooth] coordinates {(B) (C) };
\draw[black, line width=2pt] (C) arc (270:180:\r) ;
\draw [black, line width=2pt] plot [smooth] coordinates {(D)  (E)  };

\draw[black, line width=2pt] (E) arc (180:90:\r) ;
\draw [black, line width=2pt] plot [smooth] coordinates {(F) (J)  };

\end{tikzpicture}%

%% file: VisualisierungPolstellen4.tex
%
%
\begin{tikzpicture}

\tikzmath{%
\s=1;
\b=-1.5;
\d=1;
\l=2;
\r=0.52;
}

\draw[thick,->] (-4,0) -- (4,0);
\draw[thick,->] (0,-3) -- (0,3);

\tkzDefPoint(-2.5,2.5){node1}
\tkzLabelPoint[above](node1){$\lambda<0,\,\beta<0$}

\tkzDefPoint(\b,1){b1}
\tkzDefPoint(\b-\s,1){b2}
\tkzDefPoint(\b-2*\s,1){b3}
\tkzLabelPoint[above](b1){$\beta$}
\tkzLabelPoint[above](b2){$\beta-1$}
\tkzLabelPoint[above](b3){$\beta-2$}

\foreach \n in {0,...,2}
	\node at (\b-\n*\s,1)[circle,fill,inner sep=1.5pt]{};

\tkzDefPoint(\s,0){t1}
\tkzDefPoint(2*\s,0){t2}
\tkzDefPoint(3*\s,0){t3}
\tkzLabelPoint[above](t1){$1$}
\tkzLabelPoint[above](t2){$2$}
\tkzLabelPoint[above](t3){$3$}

\foreach \n in {1,...,4}
	\node at (\n*\s,0)[circle,fill,inner sep=1.5pt]{};

\tkzDefPoint(\d/\l*\s,-1){n1}
\tkzDefPoint(\d/\l*\s+2/\l*\s,-1){n2}
\tkzDefPoint(\d/\l*\s+4/\l*\s,-1){n3}
\tkzLabelPoint[above](n1){$\frac{-n}{\lambda}$}
\tkzLabelPoint[above](n2){$\frac{-n-2}{\lambda}$}
\tkzLabelPoint[above](n3){$\frac{-n-4}{\lambda}$}

\foreach \n in {1,...,4}
	\node at (-\d/\l*\s+2/\l*\s*\n,-1)[circle,fill,inner sep=1.5pt]{};

\tkzDefPoint(-2/\l*\s*1,-0.5){l2}
\tkzDefPoint(-2/\l*\s*2,-0.5){l3}
\tkzLabelPoint[below](l2){$\frac{2}{\lambda}$}
\tkzLabelPoint[below](l3){$\frac{4}{\lambda}$}

\foreach \n in {1,...,3}
	\node at (-2/\l*\s*\n,-0.5)[circle,draw,inner sep=1.5pt]{};

\coordinate (A) at (4,-1.25);
\coordinate (B) at ($(A)+ (-2,0)$) ;
\coordinate (C) at ($(A)+ (-4.5,0)+(\r,0)$);
\coordinate (D) at ($(C)+(-\r,\r)$);
\coordinate (E) at (-0.5,0);
\coordinate (F) at  ($(E)+ (\r,\r)$) ;
\coordinate (G) at  ($(F)+ (\b-0.5-\r,0)$) ;
\coordinate (H) at (\b,3*\r);
\coordinate (J) at (4,\r);

\coordinate (K) at (A)++(1,1);

\draw [black, line width=2pt, ->] plot [smooth] coordinates {(A) (B) };
\draw [black, line width=2pt] plot [smooth] coordinates {(B) (C) };
\draw[black, line width=2pt] (C) arc (270:180:\r) ;
\draw [black, line width=2pt] plot [smooth] coordinates {(D)  (E)  };

\draw[black, line width=2pt] (E) arc (180:90:\r) ;
\draw [black, line width=2pt] plot [smooth] coordinates {(F) (J)  };

\end{tikzpicture}%